# A Decomposition Result for the Haar Distribution on the Orthogonal Group

by

**Morris L. Eaton**  &  **Robb J. Muirhead**
School of Statistics    Statistical Research and Consulting Center
University of Minnesota    Pfizer Inc.

## Abstract

Let $\Gamma$ be a Haar distributed random matrix on the group $\mathcal{O}_p$ of $p \times p$ real orthogonal matrices. Partition $\Gamma$ into four blocks, $\Gamma_{11} : 1 \times 1$, $\Gamma_{12} : 1 \times (p-1)$, $\Gamma_{21} : (p-1) \times 1$ and $\Gamma_{22} : (p-1) \times (p-1)$, so

$$\Gamma = \begin{pmatrix} \Gamma_{11} & \Gamma_{12} \\ \Gamma_{21} & \Gamma_{22} \end{pmatrix}.$$

The marginal distribution of $\Gamma_{11}$ is well known. In this paper, we give the conditional distribution of $(\Gamma_{21}, \Gamma_{12})$ given $\Gamma_{11}$, and the conditional distribution of $\Gamma_{22}$ given $(\Gamma_{21}, \Gamma_{12}, \Gamma_{11})$. This conditional specification uniquely determines the Haar distribution on $\mathcal{O}_p$. The two conditional distributions involve well known probability distributions – namely, the uniform distribution on the unit sphere $\mathcal{S}_{p-1} = \{x \in R^{p-1} \mid \|x\| = 1\}$ and the Haar distribution on $\mathcal{O}_{p-2}$. Our results show how to construct the Haar distribution on $\mathcal{O}_p$ from the Haar distribution on $\mathcal{O}_{p-2}$ coupled with the uniform distribution on $\mathcal{S}_{p-1}$.

## 1. Introduction and Summary

The focus of this paper is the Haar probability distribution on the group $\mathcal{O}_p$ of $p \times p$ real orthogonal matrices. The use of this group and the Haar distribution in multivariate statistical analysis has a long history, with James (1954) and Wijsman (1957) being two important early contributions. A standard description of the Haar distribution on $\mathcal{O}_p$ is in terms of invariant differential forms – see Farrell (1985) for a systematic development and excellent history of this approach in multivariate analysis. A useful alternative is the use of random matrices, the multivariate normal distribution, and invariance properties of the objects under study. For example, see Eaton (1983) and Eaton (1989, Chapter 7). The primary technical tools used in this paper stem from the invariance considerations discussed at length in Eaton (1989).

To describe the problem under consideration in this paper, suppose the random $p \times p$ orthogonal matrix $\Gamma$ has the Haar distribution. This distribution is characterized by its invariance. To be more precise, let $\mathcal{L}(\cdot)$ denote the distribution (or probability law) of "·", where "·" can be a random variable, a random vector, a random matrix, etc. Using the $\mathcal{L}$-notation, the Haar probability distribution is characterized by

$$\mathcal{L}(\Gamma) = \mathcal{L}(g_1 \Gamma) = \mathcal{L}(\Gamma g_2)$$

for all $g_1, g_2 \in \mathcal{O}_p$. In other words, the Haar distribution is the unique invariant (right or left) probability distribution on $\mathcal{O}_p$.

In all that follows, we will assume that $\Gamma \in \mathcal{O}_p^+$, where

$$\mathcal{O}_p^+ = \left\{ h \in \mathcal{O}_p \,\middle|\, h = \begin{pmatrix} h_{11} & h_{12} \\ h_{21} & h_{22} \end{pmatrix},\, h_{11} \in (-1, 1) \right\}.$$

Note that $\mathcal{O}_p - \mathcal{O}_p^+$ is a set of Haar probability zero. In the arguments below, this set of probability zero has been removed from the sample space of $\Gamma$.

To describe the results in this paper, partition $\Gamma \in \mathcal{O}_p^+$ as

$$\Gamma = \begin{pmatrix} \Gamma_{11} & \Gamma_{12} \\ \Gamma_{21} & \Gamma_{22} \end{pmatrix}$$

where $\Gamma_{11}$ is $1 \times 1$, $\Gamma_{12}$ is $1 \times (p-1)$, $\Gamma_{21}$ is $(p-1) \times 1$ and $\Gamma_{22}$ is $(p-1) \times (p-1)$. The marginal distribution of $\Gamma_{11}$ is well known – see below following Theorem 1.1. (It is well known even if $\Gamma_{11}$ is $s \times t$ with $s + t \leq p$; see Mitra (1970), Khatri (1970) and Eaton (1989, Chapter 7).) Thus, we will proceed with $\mathcal{L}(\Gamma_{11})$ being specified. In what follows, the notation $\mathcal{L}(\cdot|*)$ is used for the conditional distribution of "·" given "*". The basic results in this paper provide a complete description of the two conditional distributions

$$\mathcal{L}(\Gamma_{21}, \Gamma_{12} | \Gamma_{11}) \tag{1.1}$$

and

$$\mathcal{L}(\Gamma_{22} | \Gamma_{21}, \Gamma_{12}, \Gamma_{11}) \tag{1.2}$$

A moment's reflection will convince the reader that knowing $\mathcal{L}(\Gamma_{11})$, (1.1) and (1.2) determines the Haar distribution on $\mathcal{O}_p$ and conversely.

Here is a rigorous specification of (1.1). Let $U_1$ and $U_2$ be independent, identically distributed (*iid*) and uniform on the unit sphere $\mathcal{S}_{p-1} = \{x \in R^{p-1} \mid \|x\| = 1\}$.

**Theorem 1.1** In the notation above and with $\Gamma_{11} \in (-1, 1)$ fixed,

$$\mathcal{L}(\Gamma_{21}, \Gamma_{12} \mid \Gamma_{11}) = \mathcal{L}\left((1-\Gamma_{11}^2)^{1/2} U_1, (1-\Gamma_{11}^2)^{1/2} U_2'\right) \tag{1.3}$$

where $U_2'$ is the transpose of $U_2$.

The above result asserts that $\mathcal{L}(\Gamma_{21}, \Gamma_{12}, \Gamma_{11})$ can be generated as follows:

(i) First, draw $\Gamma_{11}$ from the density (see Eaton (1989, Proposition 7.3))

$$f(x \mid p) = \frac{\Gamma(\frac{1}{2} p)}{\Gamma(\frac{1}{2}) \Gamma(\frac{1}{2}(p-1))} (1-x^2)^{(p-3)/2}, \qquad |x| < 1$$

(ii) Next, draw $U_1$ and $U_2$ which are *iid* uniform on $\mathcal{S}_{p-1}$. Then use (1.3) to specify the conditional distribution of $(\Gamma_{21}, \Gamma_{12})$ given $\Gamma_{11}$.

It is obvious that (i) and (ii) determine $\mathcal{L}(\Gamma_{21}, \Gamma_{12}, \Gamma_{11})$.

Our next task is to specify the conditional distribution (1.2). To this end, fix the values of $\Gamma_{11}, \Gamma_{21}, \Gamma_{12}$. and recall that $\Gamma_{11} \in (-1, 1)$. Let $h_1 \in \mathcal{O}_{p-1}$ be an orthogonal transformation satisfying

$$h_1 \varepsilon_1 = \frac{1}{\sqrt{1-\Gamma_{11}^2}} \Gamma_{21}, \tag{1.4}$$

where

$$\varepsilon_1 = \begin{pmatrix} 1 \\ 0 \\ \vdots \\ 0 \end{pmatrix} \in R^{p-1}.$$

Also, let $h_2 \in \mathcal{O}_{p-1}$ satisfy

$$h_2 \varepsilon_1 = \frac{1}{\sqrt{1-\Gamma_{11}^2}} \Gamma'_{12}. \tag{1.5}$$

That $h_1$ and $h_2$ exist and depend only on the value of $(\Gamma_{21}, \Gamma_{12})$ is demonstrated in Proposition A.2 in the Appendix. Finally, let the $(p-2) \times (p-2)$ random matrix $\Delta$ have the Haar distribution on $\mathcal{O}_{p-2}$.

**Theorem 1.2:** In the notation above with $\Gamma_{11}, \Gamma_{21}, \Gamma_{12}$ fixed, a version of the conditional distribution of $\Gamma_{22}$ given $(\Gamma_{21}, \Gamma_{12}, \Gamma_{11})$ is the distribution

$$\mathcal{L}\left( h_1 \begin{pmatrix} -\Gamma_{11} & 0 \\ 0 & \Delta \end{pmatrix} h'_2 \right) \tag{1.6}$$

where $\Delta$ is Haar distributed on $\mathcal{O}_{p-2}$. Here $h_1$ and $h_2$ are given in (1.4) and (1.5) respectively.

The proofs of both Theorem 1.1 and Theorem 1.2 rely to a large extent on some fairly well known notions from group theory and invariant measures. In Section 2, we present the underlying group action that provides the appropriate setting for our proofs. It is assumed that the reader is somewhat familiar with the standard notions of left group action, existence and uniqueness of invariant measures in the compact case, and the basic representation result given in Theorem 4.4 in Eaton (1989, Chapter 4). The proof of Theorem 1.1, given in Section 3, involves little more than the standard assertion that the uniform probability distribution on $\mathcal{S}_{p-1}$ is the unique orthogonally invariant probability measure on $\mathcal{S}_{p-1}$.

Our proof of Theorem 1.2 is somewhat more involved. It depends on a general constructive method for describing an invariant conditional distribution, given the value of an equivariant statistic. This method, which we believe is new, is presented rather abstractly in the first portion of Section 3. A direct application of the method provides a proof of Theorem 1.2.

Finally, we note here that we obtained have versions of Theorems 1.1 and 1.2 for the case where $\Gamma_{11}$ is $q \times q$, with $1 < q < p/2$. The results are available from the authors.

## 2. A Group Action on $\mathcal{O}_p$

We begin this section with a description of an invariance property of the conditional distribution $\mathcal{L}(\Gamma | \Gamma_{11})$ on $\mathcal{O}_p^+$. As above, $\Gamma$ has the Haar distribution on $\mathcal{O}_p^+$ and $\Gamma$ has

been partitioned as in Section 1. Let $\mathcal{H}$ be the compact matrix group whose $p \times p$ elements $h$ have the form

$$h = \begin{pmatrix} 1 & 0 \\ 0 & h_1 \end{pmatrix}, \quad h_1 \in \mathcal{O}_{p-1}.$$

Obviously $\mathcal{H}$ is a subgroup of $\mathcal{O}_p$ so that

$$\mathcal{L}(\Gamma) = \mathcal{L}(h\Gamma k') \tag{2.1}$$

for all $h, k \in \mathcal{H}$. The reason for the transpose on $k$ in (2.1) is so that the action of the product group $\mathcal{H} \otimes \mathcal{H}$ on $\mathcal{O}_p^+$ given by

$$\Gamma \to h\Gamma k' \tag{2.2}$$

is in fact a left action. (See Eaton (1989, pp.19-20) for the distinction between left and right actions.) The action (2.2) can be expressed in terms of the blocks of $\Gamma$ and the two $(p-1) \times (p-1)$ lower right blocks of $h$ and $k$:

$$h = \begin{pmatrix} 1 & 0 \\ 0 & h_1 \end{pmatrix}, \quad k = \begin{pmatrix} 1 & 0 \\ 0 & k_1 \end{pmatrix}.$$

The action on the blocks is

$$\begin{aligned} \Gamma_{11} &\to \Gamma_{11} \\ \Gamma_{21} &\to h_1 \Gamma_{21} \\ \Gamma_{12} &\to \Gamma_{12} k_1' \\ \Gamma_{22} &\to h_1 \Gamma_{22} k_1' \end{aligned} \tag{2.3}$$

For each $\gamma \in (-1, 1)$, let $\mathcal{X}_\gamma$ be the subset of $\mathcal{O}_p^+$ defined by

$$\mathcal{X}_\gamma = \{\Gamma \in \mathcal{O}_p^+ \mid \Gamma_{11} = \gamma\} \tag{2.4}$$

It is clear that $\mathcal{H} \otimes \mathcal{H}$ acts on $\mathcal{X}_\gamma$ for each $\gamma$, with the action being given by (2.2), or equivalently, (2.3).

Our first result implies that the action of $\mathcal{H} \otimes \mathcal{H}$ on $\mathcal{X}_\gamma$ is transitive.

**Proposition 2.1:** Consider $\psi \in \mathcal{O}_p^+$ with

$$\psi = \begin{pmatrix} \psi_{11} & \psi_{12} \\ \psi_{21} & \psi_{22} \end{pmatrix}$$

partitioned as $\Gamma$ is partitioned. Given $\psi$, define $\psi_0$ by

$$\psi_0(\psi_{11}) \equiv \psi_0 = \begin{pmatrix} \psi_{11} & \sqrt{1-\psi_{11}^2}\,\varepsilon_1' \\ \sqrt{1-\psi_{11}^2}\,\varepsilon_1 & \psi_{22}^* \end{pmatrix} \qquad (2.5)$$

where

$$\psi_{22}^* = \begin{pmatrix} -\psi_{11} & 0 \\ 0 & I_{p-2} \end{pmatrix} \qquad (p-1)\times(p-1).$$

Then there is an $h \otimes k \in \mathcal{H} \otimes \mathcal{H}$ such that

$$(h \otimes k)\psi = \psi_0,$$

where $(h \otimes k)\psi$ is specified by (2.3).

**Proof:** In the notation of (2.3),

$$(h \otimes k)\psi = \begin{pmatrix} \psi_{11} & \psi_{12}k_1' \\ h_1\psi_{21} & h_1\psi_{22}k_1' \end{pmatrix}. \qquad (2.6)$$

First pick $h_1 = \tilde{h}_1$ so that

$$\tilde{h}_1\psi_{21} = \sqrt{1-\psi_{11}^2}\,\varepsilon_1$$

and pick $k_1 = \tilde{k}_1$ so that

$$\psi_{12}\tilde{k}_1' = \sqrt{1-\psi_{11}^2}\,\varepsilon_1'.$$

Then

$$\tilde{\psi} \equiv (\tilde{h} \otimes \tilde{k})\psi = \begin{pmatrix} \psi_{11} & \sqrt{1-\psi_{11}^2}\,\varepsilon_1' \\ \sqrt{1-\psi_{11}^2}\,\varepsilon_1 & \tilde{h}_1\psi_{22}\tilde{k}_1' \end{pmatrix}.$$

The fact that $\tilde{\psi} \in \mathcal{O}_p$ now implies that $\tilde{\psi}_{22} = \tilde{h}_1\psi_{22}\tilde{k}_1'$ can be written

$$\tilde{\psi}_{22} = \begin{pmatrix} -\psi_{11} & 0 \\ 0 & \tilde{\Delta}_{22} \end{pmatrix}$$

where $\tilde{\Delta}_{22} \in \mathcal{O}_{p-2}$. Next pick $h^* \in \mathcal{H}$ to be

$$h^* = \begin{pmatrix} I_2 & 0 \\ 0 & \tilde{\Delta}'_{22} \end{pmatrix}$$

and pick $k^* = I_p \in \mathcal{H}$. Then a routine calculation shows that

$$(h^* \otimes k^*)(\tilde{h} \otimes \tilde{k})\psi = \psi_0$$

where $\psi_0$ is given in (2.5). This completes the proof.

The following three propositions are easy consequences of Proposition 2.1 coupled with standard invariance techniques. For some background material see Eaton (1989) (Section 2.3 for a discussion of maximal invariants and orbits, Section 4.1 for material on cross sections, and Theorem 4.1 for a representation result).

**Proposition 2.2:** Define $f$ on $\mathcal{O}_p^+$ by

$$f(\psi) = \psi_0(\psi_{11}) \in \mathcal{O}_p^+.$$

Then $f$ is a maximal invariant under the action of $\mathcal{H} \otimes \mathcal{H}$.

**Proof:** The invariance of $f$ is obvious. Consider $\psi$ and $\xi$ in $\mathcal{O}_p^+$ and suppose

$$f(\psi) = \psi_0(\psi_{11}) = f(\xi),$$

so that $\psi_{11} = \xi_{11}$. Use Proposition 2.1 to pick $h \otimes k$ and $\tilde{h} \otimes \tilde{k}$ so that

$$(h \otimes k)\psi = \psi_0(\psi_{11}) = (\tilde{h} \otimes \tilde{k})\xi.$$

Then $(\tilde{h} \otimes \tilde{k})^{-1}(h \otimes k)\psi = \xi$ so that $f$ is a maximal invariant. This completes the proof.

**Proposition 2.3:** The set

$$\{\psi_0(\psi_{11}) \mid \psi_{11} \in (-1, 1)\}$$

is a measurable cross section of $\mathcal{O}_p^+$.

**Proof:** This is obvious from Propositions 2.1 and 2.2.

It is now clear that the orbit of a point $\psi$ in $\mathcal{O}_p^+$ is just $\mathcal{X}_{\psi_{11}}$, and $\mathcal{H} \otimes \mathcal{H}$ acts transitively on each $\mathcal{X}_\gamma$.

Finally, let $U$ have the uniform (Haar) distribution on the compact group $\mathcal{H} \otimes \mathcal{H}$. Obviously, $U$ has the same distribution as $V_1 \otimes V_2$ where $V_1$ and $V_2$ are *iid* Haar on $\mathcal{H}$. Theorem 4.4 in Eaton (1989) immediately implies the following:

**Proposition 2.4:** For fixed $\Gamma_{11} \in (-1, 1)$,

$$\mathcal{L}(U\psi_0(\Gamma_{11})) \tag{2.7}$$

serves as a version of the conditional distribution

$$\mathcal{L}(\Gamma | \Gamma_{11}) \tag{2.8}$$

when $\Gamma$ is Haar on $\mathcal{O}_p^+$.

Note that the distribution (2.7) is on $\mathcal{O}_p^+$ but of course is concentrated on the orbit of $\psi_0(\Gamma_{11})$. The distribution (2.8) can be interpreted as a distribution on $\mathcal{X}_{\Gamma_{11}}$, or equivalently, as a distribution concentrated on the orbit of $\psi_0(\Gamma_{11})$. In all that follows, we will take (2.7) to be a version of (2.8) and will treat (2.8) as a distribution on $\mathcal{X}_{\Gamma_{11}}$. An immediate consequence of Proposition 2.4 is the following:

**Corollary 2.5:** The conditional distribution (2.8) on $\mathcal{X}_{\Gamma_{11}}$ is invariant under the action of $\mathcal{H} \otimes \mathcal{H}$ on $\mathcal{X}_{\Gamma_{11}}$; that is,

$$\mathcal{L}(\Gamma | \Gamma_{11}) = \mathcal{L}(h\Gamma k' | \Gamma_{11}) \tag{2.9}$$

for all $h, k \in \mathcal{H}$. In particular, by marginalization,

$$\mathcal{L}(\Gamma_{21}, \Gamma_{12} | \Gamma_{11}) = \mathcal{L}(h_1\Gamma_{21}, \Gamma_{12}k_1' | \Gamma_{11}) \tag{2.10}$$

for all $h, k \in \mathcal{H}$.

The proof of Theorem 1.1 is now a routine argument using standard invariance techniques. Given $\Gamma_{11}$, the vectors $\Gamma_{21}$ and $\Gamma_{12}$ satisfy $\|\Gamma_{21}\|^2 = \|\Gamma_{12}\|^2 = 1-\Gamma_{11}^2$. Let $\mathcal{Z}$ be this space of values for $(\Gamma_{21},\Gamma_{12})$. The distribution specified by (2.7) is invariant under the action of the compact group $\mathcal{D} = \mathcal{H} \otimes \mathcal{H}$. To conclude uniqueness of this distribution and the conclusion of Theorem 1.1, note that $\mathcal{Z}$ is a topological left homogeneous space under this action. Here, we are using the terminology of Nachbin (1965, p. 128) which is explained in detail at the beginning of the next section. The reader should note that $\mathcal{X}_\gamma$ is also a left homogeneous space under the action of $\mathcal{D}$ on $\mathcal{X}_\gamma$. This fact is used in the next section and in the Appendix.

## 3. An Invariant Conditional Distribution

We begin this section with a general result concerning the existence of invariant conditional distributions under rather natural invariance assumptions. This result is then used to provide a proof of Theorem 1.2.

Consider a Polish space $(\mathcal{X},\mathcal{B})$ which is acted upon topologically by a compact group $\mathcal{D}$ that is also a Polish space. Here $\mathcal{B}$ is the σ-algebra of Borel sets of $\mathcal{X}$. It is assumed that $\mathcal{X}$ is a topological left homogeneous space – see Nachbin (1965, p. 128) for a discussion of this terminology. In particular, $\mathcal{D}$ is assumed to be transitive on $\mathcal{X}$ and the map $T_x : \mathcal{D} \to \mathcal{X}$ given by $T_x(g) = gx$ is assumed to be an open mapping. Under these assumptions, there exists a unique $\mathcal{D}$-invariant probability measure $P$ defined on $\mathcal{B}$. The notation $\mathcal{L}(X) = P$ means that the random object $X \in \mathcal{X}$ has distribution $P$. Of course, $\mathcal{L}(X) = \mathcal{L}(gX)$ for all $g$ since $P$ is $\mathcal{D}$-invariant.

Next, we consider a continuous mapping $t$ from $\mathcal{X}$ onto a Polish space $(\mathcal{Y},\mathcal{C})$. It is assumed that $t$ is an *equivariant* map – that is, we assume

$$t(x_1) = t(x_2) \Rightarrow t(gx_1) = t(gx_2) \text{ for all } g \in \mathcal{D}. \qquad (3.1)$$

This assumption allows us to induce a group action on $\mathcal{Y}$. The basic idea is the following. Given $y \in \mathcal{Y}$, there is an $x \in \mathcal{X}$ such that $t(x) = y$ since $t$ is onto. Now, we simply define $gy$ to be $t(gx)$. It is assumption (3.1) that allows us to establish that this definition of $gy$, namely

$$gy = t(gx) = gt(x), \qquad (3.2)$$

is unambiguous. See Eaton (1989, Theorem 2.4 on page 32 and page 35 of Section 2.4) for details and some further discussion. In all that follows, we assume that $\mathcal{Y}$ is also a topological left homogeneous space under the action of $\mathcal{D}$.

**Remark 3.1:** As motivation for the above assumption, we note that the main application of the material to follow is for the case when $\mathcal{X}$ is $\mathcal{X}_\gamma$ given in Section 2. The mapping $t$ is given by

$$t(\Gamma) = (\gamma, \Gamma_{21}, \Gamma_{12})$$

for $\Gamma \in \mathcal{X}_\gamma$. Of course, $\mathcal{Y} = \mathcal{H} \otimes \mathcal{H}$ in this application.

The main theoretical result of this section establishes the existence of an invariant version of the conditional distribution of $X$ given $t(X) = y \in \mathcal{Y}$. More precisely, let $\mathcal{L}(X|t(X) = y)$ denote some version of the conditional distribution when $X$ has distribution $P$ above. In what follows, we will show that there is a Markov kernel, $R(\cdot|y)$ on $\mathcal{B} \times \mathcal{Y}$, that serves as a version of $\mathcal{L}(X|t(X) = y)$ and is invariant in the sense that

$$R(B|y) = R(gB|gy) \tag{3.3}$$

for all Borel sets $B \in \mathcal{B}$, for all $y \in \mathcal{Y}$ and for $g \in \mathcal{G}$. It is this invariance, when applied to the situation of Section 1, that underlies the proof of Theorem 1.2.

We now proceed with some technical details.

**Proposition 3.1:** The action of the group $\mathcal{G}$ on $\mathcal{Y}$ is transitive.

**Proof:** Consider $y_1$ and $y_2$ in $\mathcal{Y}$. We need to show that there is a $g \in \mathcal{G}$ so that $gy_1 = y_2$. Because t is an onto map, there exist $x_1$ and $x_2$ in $X$ so that $t(x_i) = y_i$ for i=1,2. But $\mathcal{G}$ is transitive on $\mathcal{X}$ by assumption, so $gx_1 = x_2$ for some $g$. Using (3.2), we have

$$y_2 = t(x_2) = t(gx_1) = gt(x_1) = gy_1$$

and the proof is complete.

**Proposition 3.2:** Let $Q = \mathcal{L}(t(X))$. Then $Q$ is an invariant probability measure on $(\mathcal{Y}, \mathcal{C})$.

**Proof:** For $C \in \mathcal{C}$ and $g \in \mathcal{G}$,

$$Q(C) = P\{t(X) \in C\} = P(t(gX) \in C) = P\{gt(X) \in C\} = P\{t(X) \in g^{-1}C\} = Q(g^{-1}C)$$

**Proposition 3.3:** For $y \in \mathcal{Y}$, let $\mathcal{X}_y = \{x \mid t(x) = y\}$. Then

$$\mathcal{X}_{gy} = g\mathcal{X}_y. \tag{3.4}$$

**Proof:** Using the equivariance of t, we have

$$\mathcal{X}_{gy} = \{x \mid t(x) = gy\} = \{x \mid g^{-1}t(x) = y\} = \{g(g^{-1}x) \mid t(g^{-1}x) = y\} = \{gu \mid t(u) = y\} = g\mathcal{X}_y.$$

This completes the proof.

Now, we proceed with the description of the transition function $R(\cdot|y)$ that is to serve as a version of $\mathcal{L}(X \mid t(X) = y)$.

(i) Fix $y_0 \in \mathcal{Y}$ and let $x_0 \in \mathcal{X}$ be any point such that $t(x_0) = y_0$.

(ii) Let $\mathcal{G}_0 \subseteq \mathcal{G}$ be the group $\mathcal{G}_0 = \{g \mid gy_0 = y_0\}$. Clearly $\mathcal{G}_0$ is a compact subgroup of $\mathcal{G}$. Let $U_0$ be the unique random element of $\mathcal{G}_0$ with the Haar (on $\mathcal{G}_0$) probability distribution.

(iii) From Proposition A.1 in the Appendix, there is a measurable map $k$ from $\mathcal{Y}$ into $\mathcal{G}$ such that $k(y)y_0 = y$ for all $y \in \mathcal{Y}$.

(iv) For $y \in \mathcal{Y}$, consider the random variable

$$Z_y = k(y)U_0 x_0 \tag{3.5}$$

and let $R(\cdot|y)$ denote the distribution of $Z_y \in \mathcal{X}$. In other words, for $B \in \mathcal{B}$,

$$R(B \mid y) = \Pr\{k(y)U_0 x_0 \in B\} = \Pr\{U_0 x_0 \in (k(y))^{-1}B\},$$

where "Pr" refers to the Haar distribution of $U_0$ on $\mathcal{G}_0$. The measurability of $k(\cdot)$ insures that $R(\cdot|\cdot)$ is a Markov kernel.

The following result establishes some basic properties of $R(\cdot|\cdot)$.

**Proposition 3.4:** The Markov kernel $R$ satisfies the following:

(i) $R(gB \mid gy) = R(B \mid y)$ (3.6)

for all $B$, $g$ and $y$.

(ii) $R(\mathcal{X}_y|y) = 1$ (3.7)

for all $y$.

**Proof:** To establish (3.6), consider fixed $g$ and $y$. Then

$$R(B|gy) = \Pr\{k(gy)U_0 x_0 \in B\} = \Pr\{g^{-1}k(gy)U_0 x_0 \in g^{-1}B\}.$$

But $k(y)$ satisfies $k(y)y_0 = y$ and $g^{-1}k(gy)y_0 = y$. Therefore $g^{-1}k(gy) = k(y)g_0$ for some $g_0 \in \mathcal{G}_0$. Thus,

$$R(B|gy) = \Pr\{k(y)g_0 U_0 x_0 \in g^{-1}B\}.$$

Since $\mathcal{L}(g_0 U_0) = \mathcal{L}(U_0)$, we conclude that

$$R(B|gy) = \Pr\{k(y)U_0 x_0 \in g^{-1}B\} = R(g^{-1}B|y).$$

Thus (3.6) holds.

For (3.7), first observe that (3.7) holds when $y = y_0$. Since $g\mathcal{X}_y = \mathcal{X}_{gy}$, (3.6) yields

$$1 = R(\mathcal{X}_{y_0}|y_0) = R(g\mathcal{X}_{y_0}|gy_0) = R(\mathcal{X}_{gy_0}|gy_0).$$

The transitivity of $\mathcal{G}$ on $\mathcal{Y}$ now gives (3.7).

Note that (3.6) implies

$$\int f(g^{-1}x)R(dx|y) = \int f(x)R(dx|g^{-1}y) \qquad (3.8)$$

for all $g \in \mathcal{G}$ and for all bounded measurable $f$. The validity of (3.8) follows from (3.6) and the standard approximation of bounded measurable functions by linear combinations of indicator functions.

**Theorem 3.1:** The Markov kernel $R$ serves as a regular version of the conditional distribution of $X$ given $t(X) = y$.

**Remark:** We are using the terminology "regular conditional distribution" in the sense defined in Section 8, Chapter 5 of Parthasarathy (1967).

**Proof:** First recall that $P$ is the unique $\mathcal{G}$-invariant probability distribution on $\mathcal{X}$. Let $K(\mathcal{X})$ denote all the bounded measurable functions on $\mathcal{X}$ and define the following integral on $K(\mathcal{X})$: For $f \in K(\mathcal{X})$, let

$$J(f) = \int_{\mathcal{Y}} \int_{\mathcal{X}} f(x) R(dx|y) Q(dy). \tag{3.9}$$

The group $\mathcal{G}$ acts on $K(\mathcal{X})$ via the action

$$(gf)(x) = f(g^{-1}x).$$

Using the invariance of $Q$ and (3.8), it is easy to show that $J(f) = J(gf)$ for all $g \in \mathcal{G}$ and $f \in K(\mathcal{X})$. Thus $J$ given by (3.9) is an invariant probability integral on $K(\mathcal{X})$. By uniqueness, we have

$$\int f(x) P(dx) = J(f) \tag{3.10}$$

for all $f \in K(\mathcal{X})$. Now, (3.10) coupled with (3.7) and Theorem 8.1 on page 147 of Parthasarathy (1967) show that $R(\cdot|y)$ is a regular conditional distribution and is unique up to sets of $Q$-measure zero. This completes the proof.

**Example 3.1 (The proof of Theorem 1.2):** Here, the general situation considered above is specialized to the case considered in Theorem 1.2. The basic idea is to identify the random variable "$Z_y$" in (3.5), since it is the distribution of $Z_y$ that provides the conditional distribution of $X$ given $t(X) = y$. To this end, we need a careful specification of the spaces involved and a clear description of the basic objects "$y_0, x_0, \mathcal{G}_0, U_0$, and $k(y)$", all of which go into the definition of $Z_y$ and its distribution.

To begin, we again let $\Gamma$ have the Haar distribution on $\mathcal{O}_p^+$ and fix the value of $\Gamma_{11}$ to be $\gamma \in (-1,1)$. From the results in Section 2, the conditional distribution of $\Gamma$ given $\Gamma_{11} = \gamma$ is concentrated on the compact set $\mathcal{X}_\gamma$ (see (2.4)) and is the unique invariant distribution under the transitive action of the compact group $\mathcal{G} = \mathcal{H} \otimes \mathcal{H}$ on $\mathcal{X}_\gamma$. For this example, the set "$\mathcal{X}$" is $\mathcal{X}_\gamma$ and is easily shown to be a left homogeneous space under the action of $\mathcal{H} \otimes \mathcal{H}$.

Next, for

$$\Gamma = \begin{pmatrix} \gamma & \Gamma_{12} \\ \Gamma_{21} & \Gamma_{22} \end{pmatrix} \in \mathcal{X}_\gamma,$$

define the map $t$ by

$$t(\Gamma) = \{\gamma, \Gamma_{21}, \Gamma_{12}\} \in \mathcal{Y}. \tag{3.11}$$

Here, $\mathcal{Y}$ is the compact space $\{\gamma\} \times \mathcal{S}_{(1)} \times \mathcal{S}_{(2)}$, where

$$\mathcal{S}_{(1)} = \{u \mid u \in R^{p-1}, \|u\|^2 = 1 - \gamma^2\} \tag{3.12}$$

and

$$\mathcal{S}_{(2)} = \{u \mid u' \in R^{p-1}, \|u\|^2 = 1 - \gamma^2\} \tag{3.13}$$

Elements of $\mathcal{S}_{(1)}$ are column vectors while elements of $\mathcal{S}_{(2)}$ are row vectors. That $t$ in (3.11) is an equivariant map is readily verified, and the action of $\mathcal{D}$ on $\mathcal{Y}$ is of course

$$\{\gamma, \Gamma_{21}, \Gamma_{12}\} \to \{\gamma, h\Gamma_{21}, \Gamma_{12}k'\} \tag{3.14}$$

for $h \otimes k = \mathcal{H} \otimes \mathcal{H} = \mathcal{D}$. The conditional distribution of $\Gamma$ given $t(\Gamma)$ is what is desired.

Now, we follow the procedure that leads to "$Z_y$". First, let

$$y_0 = \left(\gamma, \sqrt{1-\gamma^2}\,\varepsilon_1, \sqrt{1-\gamma^2}\,\varepsilon_1'\right) \tag{3.15}$$

and let

$$x_0 = \begin{pmatrix} \gamma & \sqrt{1-\gamma^2} & 0 & \cdots & 0 \\ \sqrt{1-\gamma^2} & -\gamma & 0 & \cdots & 0 \\ 0 & 0 & & & \\ \vdots & \vdots & & I_{p-2} & \\ 0 & 0 & & & \end{pmatrix} \tag{3.16}$$

It is an easy argument to show that $\mathcal{D}_0 = \{g \mid gy_0 = y_0\}$ is just $\mathcal{H}_0 \otimes \mathcal{H}_0$, where $\mathcal{H}_0 \subseteq \mathcal{H}$ is the subgroup

$$\mathcal{H}_0 = \left\{h \in \mathcal{O}_p \,\middle|\, h = \begin{pmatrix} I_2 & 0 \\ 0 & h_{22} \end{pmatrix}, h_{22} \in \mathcal{O}_{p-2}\right\}. \tag{3.17}$$

The random group element $U_0 \in \mathcal{D}_0$ with the Haar distribution is just

$$U_0 = V_{0,1} \otimes V_{0,2} \in \mathcal{H}_0 \otimes \mathcal{H}_0, \tag{3.18}$$

where $V_{0,1}$ and $V_{0,2}$ are *iid* Haar on $\mathcal{H}_0$. It is obvious that for $i = 1,2$

$$\mathcal{L}(V_{0,i}) = \mathcal{L}\left(\begin{pmatrix} I_2 & 0 \\ 0 & \Delta \end{pmatrix}\right) \tag{3.19}$$

where $\Delta$ is Haar on $\mathcal{O}_{p-2}$. We now see that with $x_0$ given by (3.16) and $U_0$ given by (3.18),

$$\mathcal{L}(U_0 x_0) = \mathcal{L}\left(\begin{pmatrix} \gamma & \sqrt{1-\gamma^2} & 0 & \cdots & 0 \\ \sqrt{1-\gamma^2} & -\gamma & 0 & \cdots & 0 \\ 0 & 0 & & & \\ \vdots & \vdots & & \Delta & \\ 0 & 0 & & & \end{pmatrix}\right) \tag{3.20}$$

where $\Delta$ is Haar on $\mathcal{O}_{p-2}$.

The next step in this proof of Theorem 1.2 is the calculation of a Borel measurable function $k(y) \in \mathcal{G}$ that satisfies

$$k(y) y_0 = y, \quad y \in \mathcal{Y} \tag{3.21}$$

for this example. With $y = \{\gamma, \Gamma_{21}, \Gamma_{12}\} \in \mathcal{Y}$, a direct application of the material in Section A.2 shows that there is a $(p-1) \times (p-1)$ orthogonal matrix $h_1$ that satisfies

$$h_1\left(\sqrt{1-\gamma^2}\varepsilon_1\right) = \Gamma_{21} \tag{3.22}$$

and $h_1$ is a Borel function of $\Gamma_{21}$. Similarly, there is a $(p-1) \times (p-1)$ orthogonal matrix $h_2$ that satisfies

$$h_2\left(\sqrt{1-\gamma^2}\varepsilon_1\right) = \Gamma'_{12} \tag{3.23}$$

and $h_2$ is a Borel function of $\Gamma_{12}$. Then, setting

$$k(y) = \begin{pmatrix} 1 & 0 \\ 0 & h_1 \end{pmatrix} \otimes \begin{pmatrix} 1 & 0 \\ 0 & h_2 \end{pmatrix}, \tag{3.24}$$

we see that (3.21) holds and $k(\cdot)$ is a Borel function of $y \in \mathcal{Y}$.

With $k(y)$ given by (3.24), the random variable $Z_y$ specified in (3.5) is

$$Z_y = \begin{pmatrix} \gamma & \Gamma_{12} \\ \Gamma_{21} & h_1 \psi_{22} h_2' \end{pmatrix} \tag{3.25}$$

where

$$\psi_{22} = \begin{pmatrix} -\gamma & 0 \\ 0 & \Delta \end{pmatrix} : \quad (p-1) \times (p-1)$$

and $\Delta$ is Haar on $\mathcal{O}_{p-2}$. By Theorem 3.1, the distribution of $Z_y$ serves as a version of the conditional distribution of $\Gamma$ given $t(\Gamma)$. This completes the proof of Theorem 1.2.

The above proof leads directly to an algorithm for generating a Haar distributed matrix $\Gamma$ on $\mathcal{O}_p$, given a Haar distributed matrix $\Delta$ on $\mathcal{O}_{p-2}$. Here is the algorithm:

1. Draw $\Gamma_{11}$ from the density $f(x|p)$ given in Section 1.

2. Next draw $U_1$ and $U_2$ iid uniform on $\mathcal{S}_{p-1}$ and let

$$\Gamma_{21} = \sqrt{1 - \Gamma_{11}^2}\, U_1,$$
$$\Gamma_{12} = \sqrt{1 - \Gamma_{11}^2}\, U_2'.$$

3. Then construct the matrices $h_1$ and $h_2$ as in (3.22) and (3.23) (by applying Proposition A.2).

Then the matrix

$$\Gamma = \begin{pmatrix} \Gamma_{11} & \Gamma_{12} \\ \Gamma_{21} & h_1 \Gamma_{22}^* h_2' \end{pmatrix}, \tag{3.26}$$

where

$$\Gamma_{22}^* = \begin{pmatrix} -\Gamma_{11} & 0 \\ 0 & \Delta \end{pmatrix}, \tag{3.27}$$

is Haar distributed on $\mathcal{O}_p$.

# Appendix

This appendix contains two technical results, the first of which establishes the existence of a measurable selector needed in the construction of the Markov kernel $R$ in Section 3.

**Section A.1:** In this section we consider a compact Hausdorff group $\mathcal{G}$ and a topological left homogeneous Polish space $\mathcal{Y}$ under the left action of $\mathcal{G}$ on $\mathcal{Y}$. Fix a point $y_0 \in \mathcal{Y}$ and let

$$\mathcal{G}_0 = \{g \mid gy_0 = y_0\}. \tag{A.1}$$

Then $\mathcal{G}_0$ is obviously a compact subgroup of $\mathcal{G}$. For each $y \in \mathcal{Y}$, let

$$\mathcal{G}_y = \{g \mid gy_0 = y\}. \tag{A.2}$$

Since $\mathcal{G}_y$ is closed, it is a compact subset of $\mathcal{G}$.

**Proposition A.1:** There exists a Borel measurable map $k$ from $\mathcal{Y}$ into $\mathcal{G}$ such that $k(y) \in \mathcal{G}_y$ for all $y \in \mathcal{Y}$. Thus $k(y)y_0 = y$ for all $y \in \mathcal{Y}$.

**Proof:** We will use the Kuratowski-Ryll-Naadzewski Theorem as stated in Aliprantis and Border (1999, p.567, Theorem 17.13). Here is a sketch of the argument. Consider the *correspondence c* defined on $\mathcal{Y}$ whose values are subsets of $\mathcal{G}$, given by $c(y) = \mathcal{G}_y$, with $\mathcal{G}_y$ defined by (A.2). (See Aliprantis and Border (1999, Chapter 16) for a discussion of correspondence.) Then $c(y)$ is a closed correspondence since each $\mathcal{G}_y$ is a closed set.

Recall that $c$ is called *weakly measurable* (see Aliprantis and Border (1999, p.558 and p.525)) if

$$V = \{y \mid c(y) \cap U \neq \phi\} \tag{A.3}$$

is a Borel subset of $\mathcal{Y}$ for each open subset $U$ of $\mathcal{G}$. But, by assumption, the map $T_{y_0}$ from $\mathcal{G}$ to $\mathcal{Y}$ defined by $T_{y_0}(g) = gy_0$ is an open mapping. Hence $T_{y_0}(U)$ is an open subset of $\mathcal{Y}$ and so is Borel.

Now, it is a routine argument to show $V = T_{y_0}(U)$ where $V$ is given by (A.3), so $V$ is open and hence Borel. Because $\mathscr{D}$ is transitive on $\mathscr{Y}$, $c(y)$ is not empty. Theorem 17.13 in Aliprantis and Border (1999) immediately implies the existence of a measurable selector $k$. (In other words, $k$ is a Borel measurable map from $\mathscr{Y}$ to $\mathscr{D}$ with $k(y) \in c(y)$ for each $y$.) This completes the proof.

**Section A.2:** Here we give an explicit formula for a symmetric orthogonal matrix that interchanges two given vectors of the same length. Consider vectors $u$ and $v$ in $R^p$ with $\|u\| = \|v\| > 0$. For $u \neq v$, set $w = u - v$ and define the $p \times p$ matrix $h$ by

$$h = \begin{cases} I_p - 2\dfrac{ww'}{w'w} & \text{if } u \neq v \\ I_p & \text{if } u = v \end{cases} \qquad (A.4)$$

**Proposition A.2:** The matrix $h$ is symmetric, orthogonal, and satisfies $hu = v$ and $hv = u$.

**Proof:** Symmetry and orthogonality are obvious. The case of $u = v$ is obvious, so assume $u \neq v$. Since $hw = -w$, $h(u - v) = -u + v$. But $h(u + v) = u + v$ because $w'(u + v) = 0$. Adding these two relations yields $hu = v$, so $hv = u$ by symmetry and orthogonality. This completes the proof.

Fix an vector $v \in R^p$ with $\|v\| > 0$ and let

$$\mathscr{S}_v = \{u \in R^p \mid \|u\| = \|v\|\}.$$

Then define $h^*$ on $\mathscr{S}_v$ to $\mathscr{O}_p$ by $h^*(u) = h$, where $h$ is given by (A.4). It is not difficult to show that the mapping $h^*$ is Borel measurable.